\documentclass[11pt]{article}
\usepackage{authblk}
\usepackage{amsmath,amssymb,amsthm}
\usepackage{listings}
\usepackage{graphicx}

\theoremstyle{definition}

\makeatletter
\renewcommand\@biblabel[1]{\relax}
\makeatother

\begin{document}

\title{On absolute central moments of Poisson distribution}
\author[1,2]{Pavel S. Ruzankin\footnote{email: ruzankin@math.nsc.ru}}
\affil[1]{\normalsize
Sobolev Institute of Mathematics, Novosibirsk, Russia}
\affil[2]{Novosibirsk State University, Novosibirsk, Russia}
\date{}
\maketitle
\vspace{-1.2cm}
\begin{abstract}
A recurrence formula for absolute central moments of Poisson distribution is suggested.

\end{abstract}

Let $X$ be a Poisson random variable with mean $m$. In this paper, we study absolute central moments ${\bf E} |X-a|^r$ about $a$ for naturals $r$.
Explicit representations for such moments may be useful in situations where
we want to test whether observations in a large sample are independent  Poisson with
given means, which is the case, for instance, for some image reconstruction techniques
in emission tomography (e.g., see Ben Bouall\`{e}gue et al., 2013 and Hebert, 1990).

An explicit representation for the mean deviation was obtained independently by Crow (1958) and Ramasubban (1958):
$${\bf E}|X-m|=2e^{-m}\frac{m^{\lfloor m \rfloor +1}}{\lfloor m \rfloor !},
$$
where $\lfloor \cdot \rfloor $ denotes
the floor function. 
%

Kendall (1943, relations (5.21) and (5.22) on p. 121) has showed that,
for all integers $r\ge 2$,
\begin{eqnarray}
{\bf E}(X-m)^r&=& m \sum_{k=0}^{r-2} {r-1 \choose k}
{\bf E}(X-m)^k,\label{kend1}\\
{\bf E}(X-m)^{r+1}&=&
rm {\bf E}(X-m)^{r-1} + m\frac{d}{dm}{\bf E}(X-m)^{r}.
\label{kend2}
\end{eqnarray}

Denote by ${}_1\! F_1$ the
confluent hypergeometric function of the first kind
$${}_1\! F_1(\alpha,\beta,z)=\sum_{n=0}^\infty
\frac{z^n}{n!} \prod_{j=0}^{n-1}\frac{ \alpha + j }{\beta + j},
$$
where $\prod_{j=0}^{-1} =1$.
Katti (1960) has derived the following representation for absolute central moments of $X$
about $a$
for all odd $r\ge1$:
\begin{equation}\label{katti}
{\bf E}|X-a|^r= - {\bf E}(X-a)^r + \frac{e^{-m}m^{\lfloor a \rfloor+1}}{
(\lfloor a \rfloor+1)! } G^{(r)}(0,0),
\end{equation}
where
$$G(\beta,t)=\exp\{ t(\lfloor a \rfloor - a+\beta +1)\}
\ {}_1\! F_1(\beta +1, \beta+\lfloor a \rfloor +2,m e^t),
$$
$G^{(r)}(\beta,t)$  is its $r$-th partial derivative with respect to $t$,
which, for $t=0$, can be computed by the recurrence formula
\begin{eqnarray*}
G^{(s+1)}(\beta,0)   &=&
(\lfloor a \rfloor - a+\beta +1) G^{(s)}(\beta, 0)
+\frac{m(\beta+1)}{
\beta + \lfloor a \rfloor +2}
G^{(s)}(\beta+1,0)
\end{eqnarray*}
for $s= 1,2,\dots$, where
$$G^{(0)}(\beta,0):=G(\beta,0)= {}_1\! F_1(\beta +1, \beta+\lfloor a \rfloor +2,m )
.
$$

The main goal of this paper is to suggest a simpler recurrence formula for
absolute central moments of $X$ about $a$.

Define the sign function
$$ {\rm sign} (y)=
\begin{cases}
-1 & \mbox{if } y\le 0, \\
1 &\mbox{if }y>0
\end{cases}
$$
and put
\begin{eqnarray*}
F(b)&:=&{\bf P} (X\le b),\\
C(r,a)&:=&{\bf E} (X-a)^r,\\
D(r,a,b)&:=&{\bf E} (X-a)^r {\rm sign} (X-b),\\
B(r,a,f)&:=&{\bf E} (X-a)^r f (X),\\
\end{eqnarray*}
where $0^0=1$ by definition, $f$ is a real-valued function.
Here $F(b)$ is the cumulative distribution function of $X$, and
$C(r,a)$ is the $r$-th central moment about~$a$.
The above definition of the sign function is not common for $y=0$, but is chosen here
for the sake of convenience, to provide the equality
\begin{equation}\label{d0}
D(0,a,b)=1-2F(b).
\end{equation}

We have
$$ {\bf E} |X-a|^r=
\begin{cases}
  C(r,a) & \mbox{if } r \mbox{ is even,} \\
  D(r,a,a) & \mbox{if } r \mbox{ is odd} .
\end{cases}
$$

\pagebreak

{\bf Theorem 1. } {\it
For all integers $r\ge 1$, reals $b\ge 0$ and $a$, and functions $f$ such that\\
${\bf E}X^r |f(X)|<\infty$, the following recurrence relations are valid:
\begin{eqnarray}
C(r,a)&=&(m-a) C(r-1,a) +
m \sum_{k=0}^{r-2} {r-1 \choose k}
C(k,a) \label{tc}\\
&=&m C(r-1,a-1) - a C(r-1,a), \label{tc2}\\
D(r,a,b)&=&(m-a) D(r-1,a,b) +
m \sum_{k=0}^{r-2} {r-1 \choose k}
D(k,a,b) \nonumber \\
&&\quad\quad
+ 2 (\lfloor b \rfloor+1 - a ) ^{r-1} e^{-m} \frac{
 m^{\lfloor b \rfloor +1} }{ \lfloor b \rfloor ! }
 \label{td}\\
&=&m D(r-1,a-1,b-1) - a D(r-1,a,b),
 \label{td2}\\
 B(r,a,f)&=&(m-a) B(r-1,a,f) +
m \sum_{k=0}^{r-2} {r-1 \choose k}
B(k,a,f) \nonumber \\
&&\quad\quad
+ m \sum_{k=0}^{r-1} {r-1 \choose k}
B(k,a,\Delta f), \label{tb}
\end{eqnarray}
where $0^0=1$ by definition,
$$\Delta f(j):=f(j+1)-f(j),$$
and all the expectations $B$ are finite.
}

\medskip

\medskip
{\bf Remark 1.}
The same approach can be used to obtain similar relations for binomial and hypergeometric
distributions.
\medskip

The following statement is a direct consequence of relations \eqref{td} and \eqref{d0}.

\medskip
{\bf Corollary 1.} {\it
\begin{eqnarray*}
{\bf E}|X-m|^3 &=& m(1-2F(m))+
2\left( (m-\lfloor m \rfloor)^2 + 2 \lfloor m \rfloor +1 \right)
e^{-m} \frac{
 m^{\lfloor m \rfloor +1} }{ \lfloor m \rfloor ! },
 \\
{\bf E}|X-m|^5& =& \left(10m^2+m\right) (1-2F(m))+
2\Big( (\lfloor m \rfloor+1-m)^4\nonumber \\
&&\quad\quad+
2m\big(2 (m-\lfloor m \rfloor)^2
+7\lfloor m \rfloor +7 -3m\big)
 \Big)
e^{-m} \frac{
 m^{\lfloor m \rfloor +1} }{ \lfloor m \rfloor ! }.
 \end{eqnarray*}
}

Notice that the cumulative distribution function $F$ can be calculated, e.g., by the built-in function
{\sf ppois()} in R or by the function {\sf scipy.stats.poisson.cdf()} in Python.

\pagebreak
{\it Proof of Theorem 1.}
First, let us prove \eqref{td}.
By Proposition~1 in Borisov and Ruzankin (2002, p. 1660),
 we have the equivalences
\begin{eqnarray*}
{\bf E}|(X-a)^r f(X)|<\infty& \Leftrightarrow & {\bf E}X^r |f(X)|<\infty
\  \Leftrightarrow \
{\bf E} |\Delta^r f(X)|<\infty
\\
& \Leftrightarrow & {\bf E}X^{r-1} |\Delta f(X)|<\infty
\  \Leftrightarrow \ {\bf E}|(X-a)^{r-1} \Delta f(X)|<\infty,
\end{eqnarray*}
where the $\Leftrightarrow$ sign means ``if and only if'',
$\Delta^r$ means applying the $\Delta$ operator $r$ times.
Hence, all the expectations $B$ in \eqref{td} are finite.
For $r\ge 1$, we have
\begin{equation*}
e^m B(r,a,f)=\sum_{j=0}^{\infty} \frac{(j-a)^r f(j) m^j}{ j!} =
\sum_{j=0}^{\infty} \sum_{k=0}^r
 {r \choose k} \frac{j^k (-a)^{r-k}  f(j) m^j}{ j!}
\end{equation*}
\begin{eqnarray}
=\ \sum_{j=0}^{\infty} \sum_{k=1}^r
  {r-1 \choose k-1}  \frac{ j^k  (-a)^{r-k}  f(j) m^j }{ j! } &+&
\sum_{j=0}^{\infty} \sum_{k=0}^{r-1}
  {r-1 \choose k} \frac{ j^k (-a)^{r-k} f(j) m^j }{ j! }
  \nonumber\\
=\ \sum_{j=0}^{\infty}
 \frac{ j (j-a)^{r-1} f(j)  m^{j} }{ j! }
&-& \sum_{j=0}^{\infty}
 \frac{a (j-a)^{r-1}  f(j)  m^{j} }{ j! }
 \nonumber\\
=\ \sum_{j=0}^{\infty}
 \frac{  (j+1-a)^{r-1} f(j+1)  m^{j+1} }{ j! }
&-&  a e^m B(r-1,a,f)
 \label{lexb}
\end{eqnarray}
\begin{eqnarray*}
&=& \sum_{j=0}^{\infty}
 \frac{  (j+1-a)^{r-1} \left(f(j) + \Delta f(j) \right)  m^{j+1} }{ j! }
 - a e^m B(r-1,a,f)
\\
&=&
m \sum_{k=0}^{r-1} {r-1 \choose k}
\sum_{j=0}^{\infty } \frac{(j-a)^k \left(f(j) + \Delta f(j) \right) m^j }{ j! }
- a e^m B(r-1,a,b)
\\
&=&m \sum_{k=0}^{r-1} {r-1 \choose k}
e^m (B(k,a,f) + B(k,a,\Delta f) )
- a e^m B(r-1,a,b),
\end{eqnarray*}
which proves \eqref{td}.

Relation \eqref{tc} follows immediately from \eqref{tb} with $f(j)\equiv 1$. Analogously,
\eqref{lexb} implies \eqref{tc2}.

\medskip
Let us now prove \eqref{td}.
By \eqref{lexb}, for $r\ge 1$, we have
\begin{eqnarray}
e^m D(r,a,b)&=&
\sum_{j=0}^{\infty}
 \frac{  (j+1-a)^{r-1} {\rm sign} (j+1-b)  m^{j+1} }{ j! }
-  a e^m D(r-1,a,b)
 \label{lex}
\end{eqnarray}
which proves \eqref{td2}.
Notice that
${\rm sign} (j+1-b) \ne {\rm sign} (j-b)$ for $j=\lfloor b \rfloor$ only. Hence,
\eqref{lex} equals
\begin{eqnarray*}
&=& \sum_{j=0}^{\infty}
 \frac{  (j+1-a)^{r-1} {\rm sign} (j-b)  m^{j+1} }{ j! }
 +
\frac{
2 (\lfloor b \rfloor+1 -a ) ^{r-1} m^{\lfloor b \rfloor +1} }{ \lfloor b \rfloor ! }
- a e^m D(r-1,a,b)
\\
&=&
m \sum_{k=0}^{r-1} {r-1 \choose k}
\sum_{j=0}^{\infty } \frac{(j-a)^k {\rm sign} (j-b) m^j }{ j! } +
\frac{
2 (\lfloor b \rfloor+1 -a ) ^{r-1} m^{\lfloor b \rfloor +1} }{ \lfloor b \rfloor ! }
- a e^m D(r-1,a,b)
\\
&=&m \sum_{k=0}^{r-1} {r-1 \choose k}
e^m D(k,a,b)
+
\frac{
2 (\lfloor b \rfloor+1 -a ) ^{r-1} m^{\lfloor b \rfloor +1} }{ \lfloor b \rfloor ! }
- a e^m D(r-1,a,b),
\end{eqnarray*}
which proves \eqref{td}.

The theorem is proved.

\end{document}